\documentclass[11pt]{amsart}
\usepackage{amssymb, latexsym, mathrsfs,   color }
\usepackage[all]{xy}
\usepackage[colorlinks=true, pdfstartview=FitV,linkcolor=blue,citecolor=blue,urlcolor=blue]{hyperref}

    \advance \topmargin by -\headheight
    \advance \topmargin by -\headsep

    \evensidemargin \oddsidemargin
\newtheorem {theorem}    {Theorem}[section]

\newtheorem {prop}[theorem]    {Proposition}

\renewcommand{\rm}{\mathrm}

\newcommand{\cal}{\mathcal}

\theoremstyle{definition}
\newtheorem{definition}[theorem]{Definition}

\numberwithin{equation}{section}

\newcommand{\la}{\langle}
\newcommand{\ra}{\rangle}
\newcommand{\res}{{\rm Res}}

\newcommand{\im}{\rm Im}

\begin{document}

\title{Linear Topological Modules over Vertex Algebras}

\author[Xuanzhong Dai,  Yongchang Zhu]{ Xuanzhong Dai,     Yongchang Zhu$^{*}$}

\address{Department of Mathematics, The Hong Kong University of Science and Technology, Clear Water Bay, Kowloon, Hong Kong}

\thanks{$^{*}$   The second author's  research is supported by Hong Kong RGC grant  16305217.}

\maketitle

\section{Introduction}

Let $\frak g $ be a semi-simple Lie algebra and let $\hat {\frak g} = \frak g [ t , t^{-1} ] \oplus {\Bbb C} K $ be the corresponding affine Kac-Moody algebra.  It is well-known that a highest weight module over $ \hat {\frak g}$ of level $k$ has the structure of a module over the level $k$ vertex algebra $V^k ( {\frak g})$ associated to $\hat {\frak g}$.
 The purpose of this note is to show that certain level $k$ representations of  $\hat {\frak g}$
 that are not in category ${\cal O}$, e.g., representations constructed in \cite{DZ}, has the structure of
  what we will call a linear topological module over vertex algebra $ V^k ( {\frak g})$.

   If $A$ is a ring, a linear topological module over $A$ is   a left $A$-module $M$ with a topology  that is invariant under translations and admits a fundamental system of neighborhoods of $0$ that consists of submodules of $M$ \cite{B}. A linear topological module over a vertex algebra $V$ is a complete linear topological module over ${\Bbb C}$ equipped with an action of vertex operators $ Y ( a , z ) = \sum a_{(n)} z^{-n-1}$, $a \in V$.
   The main ingredient in its definition is the notion of a ``field" on $M $,  which requires that all $a_{(n)}$ are continuous  and $ a_{(-n)} $ has small image for $ n $ sufficiently large.  The observation
    that $n$-th normal product of two fields is again a field makes it possible to impose Borcherds' identity.
   The precise definitions are given in Section 2.

    We expect that the vertex algebras will also play a role in the study of representations of ${\hat {\frak g}}$ that are not necessarily in category ${\cal O}$.
    As the algebra $ U (\hat {\frak g } )$ is not large enough  as a candidate for invariant differential operators on the loop groups, we need to extend it to include the coefficients of
     vertex operators from  $V^k ( {\frak g})$ so that the extended algebra still acts on useful classes of  ${\hat {\frak g}}$-modules such as ones that will appear
      in the study of automorphic forms on loop groups.   We also hope that the notion of  a linear topological module over vertex algebras
     will be useful to extend the work \cite{BKP} on loop groups  over $p$-adic fields to loop groups over ${\Bbb R}$.

  In Section 2 we will give definitions of a field on a complete linear topological vector space over ${\Bbb C}$ and linear topological modules over a vertex algebra.
  In Section 3, we show that certain function spaces over infinite dimensional spaces have the structure of linear topological modules over the vertex algebra of  $\beta-\gamma$ system.
  In Section 4, we show that the representations of affine $sl_n$ constructed in \cite{DZ} are linear topological modules over $V^1( sl_n ) $.
   In Section 5, we show the induced modules of level $k$ over $\hat {\frak g} $ under certain conditions can be completed to  linear topological modules over $V^k ( {\frak  g }) $.

  \

\section{ Definitions }

Recall that a linear topology on a vector space $M$ over ${\Bbb C}$  is a topology on $M$
 that is invariant under translations and admits a fundamental system of neighborhoods of $0$ that consists of subspaces of $M$ (e.g. \cite{B}).
 If such a topology is given on $M$, $M$ is called  a linear topological vector space.
  A sequence $\{v_n\} \subset V$ is called a Cauchy sequence if for any neighborhood $U$ of $0$, there exists an integer $N$ with the property that
\[ v_i-v_j \in U,\;\;\text{ for all } i,j\geq N.  \]
 We call $M$ complete if every Cauchy sequence converges in $ V$.

Here is a simple way to construct a linear topological space. Suppose $U_1 \supset U_2 \supset \dots $ is a sequence
 of subspaces in a vector space $M_0$ with $\cap_{i=1}^\infty U_i = \{ 0 \} $. The sequence $U_i$ defines a linear topology
  on $M_0$ such that $ U_i $'s forms a fundamental system of neighborhoods of $0$. This topology on $M_0$ is not complete in general.
   The completion $M$ of $M_0$ consists of Cauchy series
    $v_1 + v_2 + \dots $, i.e., for every $U_i $, there is $N_i$ such that $ v_j \in U_i $ for all $ j > N_i $.
     $ M$ has a fundamental system $\bar U_i $ of neighborhoods of $0$, where $\bar U_i$ consists of the
      equivalence classes of Cauchy series $\sum v_i$ in $U_i$.
    For example, $ M_0 = {\Bbb C} [ z , z^{-1} ] $, $ U_i = {\Bbb C} [ z ] z^i $, the completion $M$ of $M_0$ with respect to $\{ U_i\}$  is
     $ {\Bbb C} (( z )) $ with a fundamental system of neighborhoods $ {\Bbb C} [[ z ]] z^i $, $ i =1 , 2 , \dots $.

\begin{definition} Let $M$ be a linear topological vector space. A formal power series
\[
a(z)=\sum_{n\in \mathbb Z} a_{(n)} z^{-n-1} \in \text{End}\, V[[z,z^{-1}]]
\]
is called a field on $M$ if it satisfies two conditions:
(1)  each operator $a_{(n)}$ is continuous, i.e, for any neighborhood $U$ of $0$, there exists a neighborhood $V$ of $0$ such that $ a_{(n)} V \subset U$;
 (2) for any neighborhood $U$ of $0$, there is $N$ such that $\text{Im } a_{(-n)} \subset U$ for all $n \geq N$.
\end{definition}

For any $n \in {\Bbb Z}$, we define $n$-th normal product of two fields $a(z)=\sum a_{(n)}z^{-n-1} $ and $b(z)=\sum b_{(n)}z^{-n-1} $
 in the usual way whenever the right hand side below makes sense as a field (e.g.,\cite{K}):
\begin{equation}\label{2.1}
 a(z)_{(n)}b(z)
=\res_w(\iota_{w,z}(w-z)^n a(w)b(z)-\iota_{z,w}(w-z)^n b(z)a(w)).
\end{equation}
Note that when $M$ is complete, the right hand side of (\ref{2.1}) always makes sense as a field on $M$.  Now assume that $M$ is complete. To see the well-definedness, we write $c_{(j)}$ for the coefficient of $z^{-j-1}$ of the right hand side, and we have
\[ c_{(j)}=\sum_{i=0}^\infty (-1)^i {n \choose i} a_{(n-i)}b_{(i+j)} - \sum_{i=0}^\infty (-1)^{n+i} {n\choose i} b_{(n-i+j)}a_{(i)}.\]

If $n\geq 0$, the right hand side of the above formula is indeed a finite sum, and $c_{(j)}$ can be written as
 $c_{(j)}=\sum_{i=0}^n {n\choose i} (-1)^{n-i} [a_{(i)},b_{(j-i+n)}]$.
  Obviously $c_{(j)}$ is continuous.  For a neighborhood $U$ of $0$, we want ${\im}\,  c_{(-j)}\subset U$ when $j$ is large enough, which can be deduced from the fact that $a_{(i)}$ is continuous and ${\im} \, b_{(-j-i+n)}$ is small when $j$ is large.

If $n<0$, for $v\in M$, since $\text{Im}\,a_{(n-i)} $ and $\text{Im} \, b_{(n-i+j)}$ are small when $i$ is large, so
 $  c_{(j)} v $ is Cauchy series, thus it has a limit as $M$ is complete. This proves $ c_{(j)} $ is a linear operator on $M$.
 We can choose $N$, such that $\text{Im}\,a_{(n-i)}$ and $\text{Im}\,b_{(n-i+j)}$ are small enough when $i\geq N$, and a finite sum of continuous operators are still continuous, thus $c_{(j)}$ is continuous. To prove ${\im} \, c_{(-j)}$ is small when $j $ goes to infinity,  it suffices to show the image of $\sum_{i=0}^{N-1} (-1)^i { n \choose i} a_{(n-i)}b_{(i-j)} - \sum_{i=0}^{N-1} (-1)^{n+i} {n\choose i} b_{(n-i-j)}a_{(i)}$ is small when $j$ is large, which is true because
  $\text{Im}\,b_{(i-j)}$ and $\text{Im}\,b_{(n-i-j)}$ are small for $0\leq i\leq N-1$ and  the operator $a_{(n-i)}$ is continuous for $0\leq i \leq N-1$.

The locality of two fields on $M$ are defined as usual. It is easy to see that Dong's Lemma still holds, i.e., if fields $ a (z ) , b(z )$, and $c(z) $ are pairwise mutually local, then
 $a (z)_{(n)} b ( z ) $ and $c(z)$ are mutually local, whenever the $n$-th normal product makes sense.

A linear functional  $\phi : M \to  {\Bbb C}$ is called continuous if there exists some neighborhood $U$ of $0$, such that $\phi(U)=0$.
So a continuous linear functional is continuous in the usual sense when $\mathbb C$ is given the discrete topology.
We denote by  $ M^*$  the space of all continuous linear functionals on $M$. A continuous linear map $ f: M \to M$
induces a linear map $f^*:  M^* \to M^*$ by $ f^* (\phi ) ( v ) = \phi ( f ( v ) ) $.
 The following result maybe helpful to understand the nature of a field.

\begin{prop} \label{proposition 2.2}
Let  $a(z)=\sum_{n\in \mathbb Z} a_{(n)} z^{-n-1}$ be a field on a complete linear topological vector space $M$,  then
\[
a^\ast(z):=\sum_{n\in \mathbb Z} a_{ (-n)}^* z^{-n-1} \in \text{End}( M^* )[[z,z^{-1}]]
\]
 is a field on $M^*$ in the ordinary sense.
\end{prop}

\noindent {\it Proof:} For $\varphi \in M^*$, there exists some neighborhood $U$ of $0$ such that $\varphi (U)=0$.
Now we will show that $a^\ast_{(-n)} \varphi=0$ for sufficiently large $n$. Since $a(z)$ is a field, there exists $N\geq 0$ such that $\text{Im}\, a_{(-n)} \subset U $ for
 $n\geq N$. For  $v\in M$,
 $  ( a_{(-n)}^* \varphi) (v)= \varphi (a_{(-n)}v)=0$. This proves $ a_{(-n)}^* \varphi=0 $ for $n\geq N$,
Thus $a^\ast(z)$ is a field in the ordinary sense.\qed

\

\begin{definition} Let $V$ be a vertex algebra, a complete linear topological module over $V$ is
  a linear topological vector space $M$ equipped with a linear map
\begin{align*}
Y_M(\cdot ,z): &\;\; V \longrightarrow     {\rm End} ( M)[[z,z^{-1}]]     \\
                        &\;\; a  \longmapsto         Y_M(a,z)= \sum_{i\in \mathbb Z} a_M(i) z^{-i-1}
\end{align*}
such that the following axioms hold:
\begin{enumerate}
\item For any $a\in V$,  $Y_M(a,z)$ is a field on $M$;
 \item $Y_M(1,z)={\rm Id}$;
\item The Borcherds identity holds on $M$:
\[  Y_M ( a_{(l)} b , z ) = Y_M ( a , z )_{ ( l)} Y_M ( b , z ) \]
 which is equivalent to
 \[
 (a_{(l)}b)_{M}(n)=\sum_{i=0}^\infty (-1)^i {l \choose i} a_M(l-i)b_M(n+i)-\sum_{i=0}^\infty (-1)^{l+i} {l \choose i} b_M (n+l-i)a_M(i) ;\]
\item  \[ Y_M ( T a , z ) = \frac { d } { d z }  Y_M (  a , z ) \]
\end{enumerate}
\end{definition}

Note that the completeness of $M$ guarantees that the $l$-th normal product makes sense in the condition $(3)$.
 The notion of linear topological modules over a vertex superalgebra is defined similarly.
 And we can also define the notion of intertwining operators of linear topological modules.

An immediate example is the following. Let $M= \oplus _{i=0}^\infty M_{(i)}$ be an $\mathbb N$-graded module of a $\mathbb Z$-graded vertex algebra $V$ with the
 structure map $ a \mapsto Y_M (a , z )=\sum  a_M ( n ) z^{-n-1}$.
We let    $N_n=\oplus_{i\geq n} M_{(i)}$, it is clear that $ \cap_{n=0}^\infty N_n = \{ 0 \}$.
Now we introduce a linear topology on $M$ with $N_n$'s as a basis of neighborhoods of $0$.  The completion $\tilde M$
of  $M$ with respect to this topology  is $\tilde M =\Pi_{i=1}^\infty M_{(i)}$ with a basis of neighborhoods of $0$ as
 $\tilde N_n  = \Pi_{i\geq n } M_{(i)}$. Then it is easy to see that $a_M ( n ) $ can be extended to a continuous map
  $ a_{\tilde M} ( n) : \tilde M \to \tilde M$ and $Y_{\tilde M} (a , z ) = \sum  a_{\tilde M} ( n ) z^{-n-1} $
   is a field on $\tilde M$. It is straightforward to prove that $\tilde M$ is a linear topological module over $V$.

\

\section {linear topological modules over  $\beta -\gamma$ systems }

Let $W$ be a $2g$-dimensional real symplectic space.  The loop space $ W(( t)) = W \otimes {\Bbb R} (( t)) $
 is also a symplectic space with symplectic form given by
    \[ \la  a f (x ) , b g ( t) \ra =    \la a , b \ra \; {\rm Res} ( f ( t ) g ( t ) ).\]
 It is clear that $ W [ t^{-1} ] t^{-1} $ and $ W [[ t]] $ are Lagrangian subspaces.
  The Heisenberg group associated to the symplectic space $W(( t))$ is
  $ H = W (( t)) \times {\Bbb R} $ with multiplication
  \[    ( v_1 , c_1 ) ( v_2 , c_2 ) = ( v_1 + v_2 ,  \frac 12 \la v_1 , v_2 \ra + c_1 + c_2 ) . \]
  The function space on a Lagrangian subspace is usually a representation of $H$.
  We consider the Lagrangian subspace  $ W [ t^{-1} ] t^{-1}$.
  A function $ f :   W [ t^{-1} ] t^{-1} \to {\Bbb C}$ is called a smooth function if its restriction to
  each finite dimensional subspace is smooth. We denote the space of smooth functions
   on $ W [ t^{-1} ] t^{-1}$ by $ C^{\infty} ( W [ t^{-1} ] t^{-1} ) $.
  The Heisenberg group $H$ acts on  $ C^{\infty} ( W [ t^{-1} ] t^{-1} ) $ as follows:
  for $ v \in W [ t^{-1} ] t^{-1} $,  $v' \in W [[ t]] $, and $ c \in {\Bbb R} $,
  view them as elements in $H$ by $ ( v , 0 ) , ( v' , 0 ) $, and $ ( 0  , c ) $,
 \begin{equation} \label{3.1}
  v \cdot  f    ( x ) =  f ( x + v ) ,  \; \; v' \cdot  f ( x ) =  e^{ 2 \pi i  \la x , v' \ra } f ( x )  ,  \; \; c \cdot f ( x )
   = e^{ 2 \pi i c }   f ( x ).
 \end{equation}

The  Lie algebra $ \frak h $ of  $H$ is  $  W (( t)) \times {\Bbb R} $ with the Lie Bracket
\[     [  a f ( t ) ,  b  g ( t ) ] =   \la a , b \ra {\rm Res} ( f g )     \]
The exponential map $ {\rm exp} : {\frak h} \to H $ is
\[   ( c , a f ( t ) ) \mapsto   ( c , a f ( t ) ) .  \]
 It is easy to see that the Campbell-Baker-Hausdorff  formula holds
 \[   {\rm exp } ( x ) {\rm exp } ( y ) = {\rm exp} ( x + y + \frac 1 2  [ x , y ] ). \]

 To write down the action of $ \frak h $ on $ C^\infty ( W [ t^{-1} ] t^{-1} )$, we take
  a symplectic basis $e_1 , \dots , e_{2g} $ of $W$,
i.e., \[ \la e_{ i } , e_{ g+i } \ra = - \la e_{ g+i } , e_{ i }\ra = 1 \]
 and other pairings are $0$.  We write a vector  $ x\in W [ t^{-1} ] t^{-1} $
  as
  \[  x = \sum_{j=1}^{2g} \sum_{ n=1}^\infty   x_{ j , -n}  e_{j} t^{-n} ,\]
   almost all the coordinates $ x_{ j , -n} $ are $0$. We write a function $ f(x) \in  C^\infty ( W [ t^{-1} ] t^{-1} )$
    as
     \[  f ( x_{1, -1} , x_{2, -1}, \dots , x_{2g, -1} , x_{1, -2} , \dots ) .\]
 Then the action of $H$ is given as follows.
 For $n\geq 1$, $ e_j t^{-n} $ acts as the translation $ x_{ j , - n } \mapsto x_{j, -n} +1 $.
 For $ n\geq 0 $,  $ 1\leq j \leq g $, $ e_j t^{n} $ acts as the multiplication operator
  $ e^{ - 2 \pi i x_{ g+j , - n -1} } $ and $ e_{g+j} t^{n} $ acts as multiplication operator
  $ e^{  2 \pi i x_{ j , - n -1} } $.  The central element $c\in H $ acts as the multiplication by the
   scalar $ e^{ 2 \pi i c }$.
 It follows that the Heisenberg Lie algebra $\frak h$ acts as follows.
For $n\geq 1 $,
\begin{eqnarray}\label{3.2}
&&   e_j t^{-n}  =   \frac { \partial} {\partial   x_{j,-n } }  ,  1\leq j \leq 2g  \\
&&   e_j t^{n-1}  =     - 2 \pi i x_{g+j , - n},    1 \leq  j \leq g   \nonumber \\
&&  e_{g+j} t^{n-1} =  2 \pi i x_{j , - n}  \nonumber
\end{eqnarray}
 and the center $ c  $ acts as $ 2 \pi i c   \, {\rm Id} $.

For $ n \geq 1 $ and $ 1\leq i \leq g $, we write
 \begin{align*}
  \beta_{i}(n-1):=& e_{i+g} t^{-n}=\frac {\partial} { \partial_{x_{i+g,-n}}} \\
    \beta_i(-n):= & -e_{i+g} t^{n-1} = -2\pi i x_{i,-n} \\
   \gamma_{i} (n-1):=&  e_{i} t^{-n}=\frac {\partial} {\partial_{x_{i,-n}}}\\
  \gamma_{i}(-n):= &- e_i t^{n-1} =2\pi i x_{i+g, -n}
    \end{align*}
Then we have the commutation relation
\begin{equation}
  [ \beta_{k} ( m ) , \beta_j ( n ) ] = 0 , \; \; [ \gamma_{k} ( m ) , \gamma_j ( n ) ] = 0, \; \;
   [ \beta_{k} ( m ) , \gamma_j ( n ) ] = \delta_{ k , j } \delta_{ m+ n + 1 , 0 } 2 \pi i
\end{equation}

\

We give a linear topology on $ C^\infty ( W [ t^{-1}] t^{-1} )$.
Let
 \[  R_n : C^\infty ( W [ t^{-1}] t^{-1} ) \to C^\infty ( W t^{-1} + \dots + Wt^{-n}  )\]
 be the restriction map, $ U_n :={\rm Ker} ( R_n )$, it is clear that $ U_n \supset U_{n+1}$ and
  $ \cap_{n=1}^\infty U_n =0 $. The $U_n$'s define a linear topology on $ C^\infty ( W [ t^{-1}] t^{-1} )$
  which is easily seen to be complete.
 The generating series
 \begin{equation}\label{3a}
  \beta_i ( z ) = \sum_{n=-\infty}^\infty \beta_i ( n) z^{-n-1 } , \; \;
  \gamma_i ( z ) = \sum_{n=-\infty}^\infty \gamma_i ( n) z^{-n-1 },   \;  1\leq  i \leq g
  \end{equation}
  are fields on the linear topological space $ C^\infty ( W [ t^{-1}] t^{-1} )$
  and we have
\begin{prop}
  $ C^\infty ( W [ t^{-1}] t^{-1} )$ is a linear topological module over the vertex algebra of $\beta-\gamma $  systems
under the action of (\ref{3a}).
  \end{prop}

  Here the vertex algebra of $\beta-\gamma $  systems means the vertex algebra with
  the following operator product expansion:
\begin{equation}
    \beta_k ( z ) \beta_j ( w )  \sim 0 , \;  \gamma_k ( z ) \gamma_j ( w )  \sim 0, \; \beta_k ( z ) \gamma_j ( w )
     \sim \frac { 2 \pi i } { z - w } \delta_{kj} .
  \end{equation}
Our $\beta-\gamma$ system here differs from the usual one by the factor $2\pi i$ in the third relation, which is caused
 by $2\pi i $ in the second and third equations in (\ref{3.1}). It would be the same as the usual one if we replace $2\pi i  $
  in (\ref{3.1}) by $ 1 $.   We choose to use $2\pi i  $ instead of $ 1 $ in (\ref{3.1}) for the following reason.
  In $ C^\infty ( W [ t^{-1}] t^{-1} )$, there is a subspace ${\cal S} ( W [ t^{-1}] t^{-1} )$
   which consists of Schwartz functions. Recall that $ f \in C^\infty ( W [ t^{-1}] t^{-1} )$ is a Schwartz function
    if the restriction of $f$ to any finite dimensional subspace in $W [ t^{-1}] t^{-1}$ is a Schwartz function. The definition of a Schwartz function on a finite dimensional space and infinite dimensional space is stated later in Section 4.
 The subspace ${\cal S} ( W [ t^{-1}] t^{-1} )$ is stable under the action (\ref{3.1}), however it is not
  without the factor $2 \pi i $. It is clear that ${\cal S} ( W [ t^{-1}] t^{-1} )$ is
   stable under the Lie algebra action of $\frak h$. The space  ${\cal S} ( W [ t^{-1}] t^{-1} )$
    is a complete linear topological space with a system of neighborhood $ {\cal S} ( W [ t^{-1}] t^{-1} )\cap U_n $, $ n= 1 , 2 \dots $.
     Notice that $ {\cal S} ( W [ t^{-1}] t^{-1} )\cap U_n$  is the kernel of the restriction map
     $ {\cal S} ( W [ t^{-1}] t^{-1} )\to {\cal S} ( W t^{-1} + \dots + W t^{-n} ) $.

  \begin{prop}
  $ {\cal S} ( W [ t^{-1}] t^{-1} )$ is a linear topological module over the vertex algebra of $\beta-\gamma $  systems
under the action of (\ref{3a}).
  \end{prop}

  It is known that the quadratic expressions $ : \beta_i ( z) \beta_j ( z):$, $ : \gamma_i ( z) \gamma_j ( z):$,
   $ : \beta_i ( z) \gamma_j ( z):$ generates a vertex subalgebra $V$ of $\beta - \gamma $ system that is isomorphic to
  a  vertex algebra of the affine Lie algebra of $sp_{2g}$ (e.g. \cite{FF}). Therefore
     $ C^\infty ( W [ t^{-1}] t^{-1} )$ and $ {\cal S} ( W [ t^{-1}] t^{-1} )$ are linear topological modules over the vertex algebra
     $V$. The result in \cite{Z} implies that the representation $ {\cal S} ( W [ t^{-1}] t^{-1} )$
 can be integrated to a representation of the loop group of $ SP_{2n} ( {\Bbb R} )$,  while $ C^\infty ( W [ t^{-1}] t^{-1} )$ can not be integrated.

\

\section{An Example for Affine $sl_n $}

  We denote by  $E_{uv}$  the $n\times n $ matrix with $(u,v)$-entry $1$ and other entries $0$.
  So $ E_{uv} $ ($u\ne v$) and $ E_{i , i} - E_{i+1 , i+1 }$ form a basis of the Lie algebra $sl_n$.
  Recall that a smooth function $f$ on $ {\Bbb R}^N$ is a Schwartz function if  partial derivatives $\partial^I f $ of any  order
   are rapidly decaying in the sense that $ | \partial^I f ( x ) | \leq C_{I, m} ( 1 + | x|^2 )^{-m }$ for all $m > 0$.
   We call a function $ f :  {\Bbb R}^n [ t^{-1} ] t^{-1} \to {\Bbb C}$ a Schwartz function if the restriction of $f$ to
    every finite dimensional subspace in  ${\Bbb R} [ t^{-1} ] t^{-1} $ is a Schwartz function.
  Let $ {\cal S} ( {\Bbb R}^n [ t^{-1} ] t^{-1} ) $ denote the space of all Schwartz functions on ${\Bbb R} [ t^{-1} ] t^{-1} $.
  Similar to the Schwartz spaces in Section 3, $ {\cal S} ( {\Bbb R} [ t^{-1} ] t^{-1} ) $ is a complete
   linear topological space with a system of neighborhoods of $0$ formed by the kernels of restriction
    maps $ {\cal S} ( {\Bbb R}^n  [ t^{-1} ] t^{-1} ) \to {\cal S} ( {\Bbb R}^n t^{-1} + \dots + {\Bbb R}^n t^{-m} ) $.
    We write a vector $ x \in {\Bbb R}^n  [ t^{-1} ] t^{-1}$ as
    \[  x = \left( \begin{matrix} x^1_{-1} \\ x^2_{-1} \\ \dots \\ x^n_{-1} \end{matrix} \right) t^{-1}
        + \dots  + \left( \begin{matrix} x^1_{-m} \\ x^2_{-m} \\ \dots \\ x^n_{-m} \end{matrix} \right) t^{-m} +\dots .
        \]
    and write $ f (x ) \in {\cal S} ( {\Bbb R}^n  [ t^{-1} ] t^{-1} ) $ as
    \[ f ( x^1_{-1} , x^2_{-1} , \dots , x^n_{-1} , \dots ,  x^1_{-m} , x^2_{-m} , \dots , x^n_{-m}, \dots ) .\]
    The space ${\cal S} ( {\Bbb R}^n  [ t^{-1} ] t^{-1} ) $ is closed under the multiplication operators $ x^s_{-i}$
   and differential operators $\partial_{ x^s_{-i} } $. We introduce the operator
   $\pi_t :   {\cal S} ( {\Bbb R}^n  [ t^{-1} ] t^{-1} ) \to {\cal S} ( {\Bbb R}^n  [ t^{-1} ] t^{-1} )$ by
   \begin{equation}
      \pi_t f( x^1_{-1} , x^2_{-1} , \dots  ) =  \int_{ {\Bbb R}^n } f ( y_1, \dots , y_n ,   x^1_{-1} , x^2_{-2} , \dots ) dy_1 \cdots d y_n.
  \end{equation}
  For example, for $  \phi_{c}(x)= e^{ - \pi c \sum_{ j=1}^n \sum_{k=1}^\infty  {x_{  - k }^j }^2 }$, where $ c > 0 $,
     \begin{equation}\label{4.2}
  \pi_t \phi_c ( x ) = \lambda_c \phi_c ( x ),
  \end{equation}
   where
   \[ \lambda_c = \int  e^{ - \pi c (x_1^2 + \cdots + x_n^2 ) } dx_1 \cdots dx_n = c^{-\frac  n 2 } . \]
  We have the following relations for the above operators on  $ {\cal S} ( {\Bbb R}^n  [ t^{-1} ] t^{-1} )$
  \begin{equation}
    x_{-i}^s  \pi_t = \pi_t  x_{-i-1}^s, \; \; \; \partial_{x_{-i}^s}  \pi_t = \pi_t \partial_{ x_{-i-1}^s} , \; \; \;  (\pi_t)^N  \partial_{x_{-i}^s}= 0  \; \; {\rm for}
   \; 1\leq s \leq N .
\end{equation}
  It is simple to verify the first two relations. The third relation follows from  the fact that
 \[    \int_{ \Bbb R} \partial_{x } f ( x ) d x = 0 \]
for a Schwartz function $ f \in {\cal S} ( {\Bbb R} )$.

  For each scalar $\lambda $, we introduce the eigenspace of the operator $\pi_t$:
  \begin{equation}
   {\cal S} ( {\Bbb R}^n  [ t^{-1} ] t^{-1} )_\lambda = \{    f \in {\cal S} ( {\Bbb R}^n  [ t^{-1} ] t^{-1} ) \; | \; \pi_t f = \lambda \, f \}.
  \end{equation}
  By (\ref{4.2}),  $\phi_c \in   {\cal S} ( {\Bbb R}^n  [ t^{-1} ] t^{-1} )_\lambda $ when $ \lambda_c = \lambda $.
  In \cite{DZ}, we proved the following

 \begin{prop}
  The space $ {\cal S} ( {\Bbb R}^n  [ t^{-1} ] t^{-1} )_\lambda $ is a representation of
   the affine Lie algebra $\hat {sl_n} = sl_n [ t, t^{-1} ] \oplus {\Bbb C} K $ of level $1$ under the following action,
  for $j \geq 0 $,
\begin{eqnarray}
  && \label{4.5}  \pi (  E_{uv} t^j )  =  - \sum_{ i = 1}^{\infty}   x^v_{ -i - j }   \partial_ {x^u_{ - i } } \\
  &&  \label{4.6}  \pi(E_{uv}t^{-j})=-\lambda^{-j}\pi_t^j\sum_{i=1}^\infty x_{-i}^v \partial_{x_{-i-j}^u},\\
  && \label{4.7}  \pi ((E_{uu}-E_{vv})t^{-j})=-\lambda^{-j}\pi_t^j\sum_{i=1}^\infty (x_{-i}^u \partial_{x_{-i-j}^u} -x_{-i}^v \partial _{x_{-i-j}^v })
   \end{eqnarray}
 \end{prop}

In \cite{DZ}, we  construct an action of the loop group of $ SL_n ( {\Bbb R} )$ on  $ {\cal S} ( {\Bbb R}^n  [ t^{-1} ] t^{-1} )_\lambda $
  as Hecke operators, which appears more natural than the Lie algebra action in Proposition 4.1.
 The space $ {\cal S} ( {\Bbb R}^n  [ t^{-1} ] t^{-1} )_\lambda $ has a complete linear topology
  with the kernels of restriction maps $  {\cal S} ( {\Bbb R}^n  [ t^{-1} ] t^{-1} )_\lambda \to {\cal S} (  {\Bbb R}^n t^{-1} + \dots +   {\Bbb R}^n t^{-m} ) $
   as a fundamental system of neighborhoods of $0$. It is not hard to see that operators
   in (\ref{4.5})-(\ref{4.7}) are continuous operators and the generating series $ a ( z ) $ for $ a \in sl_n $
 \begin{equation}
  a ( z ) = \sum_{ k=-\infty}^\infty \pi ( a t^k ) z^{-k-1}
\end{equation}
is a field on the complete linear topological space ${\cal S} ( {\Bbb R}^n [ t^{-1} ] t^{-1} )_\lambda $.

\begin{prop}
The linear topological space  ${\cal S} ( {\Bbb R}^n [ t^{-1} ] t^{-1} )_\lambda $ is
a linear topological module over $ V^1 (\hat { sl_n} ) $.
\end{prop}

\

\section{Induced Modules as Linear Topological Modules over Vertex Algebras}

Let $ M$ be any $\hat {\frak g}$-module of level $k$, and for each $ n > 0$, let
 $U_n := {\frak g} [ t^{-1} ] t^{-n}  M $. It is clear that $ U_n \supset U_{n+1} $.
   We assume  $\cap_{n=1}^\infty U_n =\{ 0 \}$ (if the assumption does not hold, we
    replace $M$ by the quotient module $M/ \cap_{n=1}^\infty U_n$).
 The subspaces  $\{ U_n\}$  define a linear topology on $M$,  and we denote by  $ \tilde M $ the completion of $M$.
   Since $  [ a t^k ,  {\frak g} [ t^{-1} ] t^{-n} ] \subset    {\frak g} [ t^{-1} ] t^{-n+k } $ for $ -n + k < 0 $,
    $\hat {\frak g}$-module structure on $M$ extends to $\tilde M $, and $ at^n \in \hat {\frak g}$ are continuous operators on $\tilde M$.
   It is clear that for $a \in {\frak g}$, $a ( z ) = \sum_{n=-\infty}^\infty at^n z^{-n-1} $ are fields on $\tilde M$ and
      $\tilde M$ is a linear topological module over
      $V^k (\frak g)$.

  The $\hat {\frak g}$-modules that are not in category ${\cal O}$ on which the Cartan subalgebra acts semisimply
   can be constructed using parabolic induction (see \cite{FK} and references therein).  More generally, for a subalgebra $ {\frak s} \subset \hat {\frak g}$,
    assume ${\frak s }$ contains the central element $K$ and there is a Lie algebra homomorphism $ \pi : {\frak s} \to {\Bbb C} K $ with $ \pi ( K ) = K $.
    For a complex number $k$, we have a one dimensional  ${\frak s}$-module ${\Bbb C}_k $ on which $ K$ acts as $k$ and $a\in  {\frak s}$ acts as $\pi( a ) $.
 Let $M$ in the previous paragraph be the induced module $ U ( \hat {\frak g} ) \otimes_{ U ( {\frak s} )}  {\Bbb C}_k $.
  Under the assumption that
   $ {\frak s} \cap  {\frak g} [ t^{-1} ] t^{-n} = \{ 0 \} $ for some $ n > 0 $, the condition $ \cap_{n=1}^\infty U_n =0$ is satisfied.
    Therefore the completion $\tilde M$ with respect to $\{ U_n\} $ is a linear topological module over $ V^k ( {\frak g } ) $.

    For a case related to conformal blocks on curves, we consider the completed affine Lie algebra
    \[ \tilde {\frak g } = {\frak g} \otimes {\Bbb C} (( t)) \oplus {\Bbb C} K \]
   with the Lie bracket
   \[  [ a  f (t) , b g ( t ) ] = [ a , b ] f ( t) g ( t ) +  ( a , b )  \res ( g d f ) K . \]
     For a smooth complete curve $C$ over ${\Bbb C}$ with a marked point $p \in C$ and  a local parameter $t$ at $p$, i.e., $t$ is a meromorphic function on $C$
 with $v_p ( t ) =1 $. With this choice of $t$, the fraction field of the completed local ring $\widehat {\cal O}_t $ is identified with $ {\Bbb C} (( t )) $.
Let $ A $ be the space of meromorphic functions on $ C$ that is regular on  $C-\{ p\}$,
 so $A$ is a subalgebra $\frak s$ of  $ {\Bbb C} (( t )) $.
The space $\frak g \otimes {\frak s } \oplus {\Bbb C} K $ is a subalgebra of $\tilde {\frak g} $ and the projection map
  $\frak g \otimes {\frak s } \oplus {\Bbb C} K  \to {\Bbb C} K $ is a Lie algebra homomorphism.
  If the genus of $C$ is not $0$, then ${\frak s} \cap   {\frak g } [ t^{-1} ] t^{-1} = \{ 0 \} $.
   Let  $ M = U (  \tilde {\frak g } )\otimes_{ U ( {\frak s } ) } {\Bbb C}_k $ be the induced module and $\tilde M $ be its completion.
   By the above discussion, $\tilde M $ is a linear topological module over $ V^k ( {\frak g } )$.

  A similar construction works if $\hat {\frak g}$ is replaced by the Heisenberg algebra $\frak h $ in Section 3. The linear topological modules over
   $\beta-\gamma$ systems constructed by induction for certain subalgebra ${\frak s} $ are related to the construction in Section 3. Consider
     \begin{equation} \label{5.1}
      f (  x) =  e^{ -\frac  1 2  \sum_{j=1}^{2g} \sum_{ n=1}^\infty x_{ j , -n }^2 } \in C^\infty ( W [ t^{-1} ] t^{-1} ) .\end{equation}
 It is easy to see that
 $ f ( x ) $ is killed by the operators
 \begin{equation} \label{5.2}
  x_{ j , - n } - \partial_{ x_{ j , - n } } , \; \; 1\leq  j \leq 2g , \;  n \in {\Bbb Z}_{>0}.
 \end{equation}
  Let ${\frak s} $ be the Lie subalgebra of ${\frak h }$ spanned by the operators in (\ref{5.2}) and the center.
  Then the completion of the induced module from ${\frak s}$ is isomorphic to the submodule generated by $ f ( x ) $ in (\ref{5.1}) as
   a linear topological module over the $\beta - \gamma $ system.

\end{document}